\documentclass{article}

\usepackage{amsmath}
\usepackage{amsfonts}
\usepackage{amssymb}
\usepackage{enumerate} 

\newtheorem{theorem}{Theorem}

\newtheorem{corollary}[theorem]{Corollary}

\newtheorem{lemma}[theorem]{Lemma}

\begin{document}

\title{Products of Three $k-$Generalized Lucas Numbers as Repdigits}
\author{Alaa ALTASSAN and Murat ALAN\\
Department of Mathematics, King Abdulaziz University,\\
21589, Jeddah, Saudi Arabia, aaltassan@kau.edu.sa\\
Department of Mathematics, Yildiz Technical University,\\
34210, Istanbul, Turkey, alan@yildiz.edu.tr
}

\date{}

\maketitle

\begin{abstract}
Let $ k \geq 2 $ and let  $ ( L_{n}^{(k)} )_{n \geq 2-k} $ be the  $k-$generalized Lucas sequence with certain initial $ k $ terms  and each term afterward is the sum of the $ k $ preceding terms. In this paper, we find all repdigits which are products of arbitrary three terms of $k-$generalized Lucas sequences. Thus, we find all non negative integer solutions of Diophantine equation 
$L_n^{(k)}L_m^{(k)}L_l^{(k)} =a \left( \dfrac{10^{d}-1}{9} \right)$  where  $n\geq m \geq l \geq 0$ and $1 \leq a \leq 9.$
\end{abstract}

\section{Introduction}
Let $ k \geq 2 $ be an integer. The $ k-$generalized Lucas sequence or, for simplicity, the $k-$Lucas sequence is defined by the recurrence relation as
$$ L_{n}^{(k)}=L_{n-1}^{(k)} + \cdots + L_{n-k}^{(k)} \quad \text{for all} \quad n \geq 2 $$
with the initial condition $ L_{0}^{(k)}=2$, $ L_{1}^{(k)}=1$  for all $ k \geq 2 $ and $ L_{2-k}^{(k)}= \cdots =L_{-1}^{(k)} =0 $ for all  $k \geq 3.$ For $ k=2, $ this sequence is the classical Lucas sequence, and in this case we omit the superscript $ (k) $ in the notation.

Recall that, a positive integer is called a repdigit if it has only one distinct digit in its decimal expansion. These are the numbers of the form $ a {(10^{d}-1)/9} $ for some $ d \geq 1 $ and $ 1 \leq a \leq 9 .$ In recent years, many authors have worked on problems involving relations between terms of some binary recurrence sequences and repdigits \cite{Dam,Dam21,EKZ21,Luca2012,NLT181}. Some authors extended these problems to the case involving order $ k $ generalization of these binary recurrence sequences \cite{Alahmadi,BL13,BL,14,Coufal,Marques15,24,Raya}. In fact, in \cite{Luca2000}, Luca found all repdigits in Lucas sequence whereas in \cite{16}, the authors extend this result to the $k-$Lucas sequences. Recently, in \cite{EK}, the authors determined all repdigits which are products of two Lucas numbers. In this paper, we continue to search in this line and extend some previous works, by examining the product of arbitrary three $ k-$Lucas numbers and repdigits.  More precisely, we find all repdigits which are product of three $ k-$Lucas numbers. Our main result is stated in the following theorem.
\begin{theorem}
\label{main1}
If the Diophantine equation
\begin{equation}
L_n^{(k)}L_m^{(k)}L_l^{(k)} =a \left( \dfrac{10^{d}-1}{9} \right), \quad  d \geq  2 \quad n\geq m \geq l \geq 0 \quad \text{and} \quad 1 \leq a \leq 9
\label{LnLmLlR}
\end{equation}
has solutions in non-negative integers, then
\begin{align*}
(k,n,m,l,a,d) \in & \{ (2,5,0,0,4,2), (2,5,1,0,2,2), (2,5,1,1,1,2),(2,5,2,0,6,2), \\
& (2,5,2,1,3,2),(2,5,2,2,9,2),(2,5,3,0,8,2),(2,5,3,1,4,2), \\
& (2,5,4,1,7,2), (4,5,0,0,8,2),(4,5,1,0,4,2),(4,5,1,1,2,2),\\
& (4,5,2,1,6,2) \}.
\end{align*}
\end{theorem}

\begin{corollary}
If the product of any two $k-$Lucas numbers is a repdigit with at least two digits then this product is one of the following up to order of the factors:
\begin{align*}
L_{5}^{(2)} L_{0}^{(2)}=11 \cdot 2=22, \quad & \quad L_{5}^{(2)} L_{1}^{(2)}=11 \cdot 1=11, \quad & \quad
L_{5}^{(2)} L_{2}^{(2)}=11 \cdot 3=33, \\
L_{5}^{(2)} L_{3}^{(2)}=11 \cdot 4=44, \quad & \quad L_{5}^{(2)} L_{4}^{(2)}=11 \cdot 7=77, \quad & \quad L_{5}^{(4)} L_{0}^{(4)}=22 \cdot 2=44, \\
L_{5}^{(4)} L_{1}^{(4)}=22 \cdot 1=22, \quad & \quad L_{5}^{(4)} L_{2}^{(4)}=22 \cdot 3=66. \\  
\end{align*}
\end{corollary}
In particular, we state the following result for  the  Lucas sequence.
\begin{corollary}
The only repdigits that can be written as products of three Lucas numbers are 2, 1, 3, 4, 7, 11, 22, 33, 44, 66, 77, 88 and 99 . 
\end{corollary}

The proof depends on some techniques on lower bounds for linear forms in logarithms of algebraic numbers. Also, we use a version of reduction algorithm due to Dujella and Peth\H o  \cite{DP}, which is in fact originally introduced by Baker and Davenport in \cite{Baker-Davenport}.

\section{The Tools}
Let $\theta$ be an algebraic number, and let
\[
c_0x^d+c_1x^{d-1}+\cdots+c_d=c_0\prod_{i=1}^{d}(x-\theta^{(i)})
\]
be its minimal polynomial over $ \mathbb{Z} ,$ with degree $d,$ where the $c_i$'s are relatively prime integers with $c_0>0,$ and the $\theta^{(i)}$'s are the conjugates of $\theta$.

The logarithmic height of $\theta$ is defined by
$$
h(\theta)=\dfrac{1}{d}\left(\log c_0+\sum_{i=1}^{d}\log\left(\max\{ \lvert \theta^{(i)} \rvert ,1 \}\right)\right).
$$
In particular, the logarithmic height of a rational number $\theta= r/s,$ where $ r $ and $ s $ are relatively prime integers and $s>0,$ is  $h(r/s)=\log \max \{\lvert r \rvert,s\}$. The following properties are very useful in calculation of a logarithmic height :
\begin{itemize}
\item[$ \bullet $] $h(\theta_1 \pm \theta_2)\leq h(\theta_1) + h(\theta_2)+\log 2$.
\item[$ \bullet $] $h(\theta_1 \theta_2^{\pm 1})\leq h(\theta_1)+h(\theta_2)$.
\item[$ \bullet $] $ h(\theta^{s})=\lvert s \rvert h(\theta),$ $ s \in \mathbb{Z} $.
\end{itemize}

\begin{theorem}[Matveev's Theorem, \cite{Matveev}]
\label{Matveev}
Assume that $\alpha_1, \ldots, \alpha_t$ are positive real algebraic numbers in a real algebraic number field $\mathbb{K}$ of degree $ d_\mathbb{K} ,$ and let $b_1, \ldots, b_t$ be rational integers, such that 
\[
\Lambda:=\alpha_1^{b_1}\cdots\alpha_t^{b_t}-1,
\]
is not zero. Then
$$
\lvert \Lambda\rvert >\exp\left( C(t) d_\mathbb{K}^2(1+\log d_\mathbb{K})(1+\log B)A_1\cdots A_t\right),
$$
where
\[
C(t):= -1.4\cdot 30^{t+3}\cdot t^{4.5}  \quad \text{,} \quad    B\geq \max\{\lvert b_1 \rvert,\ldots,\lvert b_t \rvert \},
\]
and
\[
A_i\geq \max\{d_\mathbb{K}h(\alpha_i),\lvert \log \alpha_i \rvert, 0.16\}, \quad \text{for all} \quad  i=1,\ldots,t.
\] 
\end{theorem}

For a real number $ \theta, $ we put $ \lvert \lvert \theta \rvert \rvert =\min\{ \lvert \theta -n \rvert :  n \in \mathbb{Z} \}, $ which represents the distance from $ \theta $ to the nearest integer. Now, we cite the following lemma which we will be used to reduce some upper bounds on the variables.
\begin{lemma} \label{reduction}\cite[Lemma 1]{13}      
Let $M$ be a positive integer, and let $p/q$ be a convergent of the continued fraction of the irrational number $\tau$ such that $q>6M.$ Let $A, B$ and $\mu$ be some real numbers with $A>0$ and $B>1$. If $\epsilon:= \lvert \lvert \mu q \rvert \rvert -M \lvert \lvert \tau q \rvert \rvert >0$, then there is no solution to the inequality
\[
0< \lvert u\tau-v+\mu \rvert <AB^{-w},
\]
in positive integers $u,v$ and $w$ with
\[
u\leq M \quad\text{and}\quad w\geq \frac{\log(Aq/\epsilon)}{\log B}.
\]
\end{lemma}

\section{Properties of $k-$Lucas Numbers }

The characteristic polynomial of the $ k-$generalized Lucas sequence is
$$\Psi_k(x)=x^k-x^{k-1}- \cdots -x-1,$$ 
which is an irreducible polynomial over $ \mathbb{Q}[x] .$ The polynomial $ \Psi_k(x) $ has exactly one real distinguished root $ \alpha(k) $ outside the unit circle \cite{22,23,26}.  The other roots of $ \Psi_k(x) $ are strictly inside the unit circle \cite{23}. This root $ \alpha(k) ,$ say for simplicity $ \alpha, $ is located in the interval
$$ 2(1-2^{-k})< \alpha < 2  \quad \text{for all} \quad k \geq 2. $$
Let
\begin{equation*}
f_k (x) =  \dfrac{x - 1}{2 + (k+1)(x - 2)}.
\end{equation*}
It is known that the inequalities
\begin{equation}\label{fkprop}
1/2 <f_k (\alpha) < 3/4  \quad \text{and} \quad \left \lvert f_k( \alpha_{i} ) \right\rvert <1, \quad 2 \leq i \leq k
\end{equation}
are hold,  where $ \alpha:=\alpha_{1}, \cdots , \alpha_{k}  $ are all the roots of $ \Psi_k(x) $  
\cite[Lemma 2]{13}. In particular, we deduce that $  f_k(\alpha) $ is not an algebraic integer. In the same lemma, it is also proved that
\begin{equation}\label{hfk}
h( f_k (\alpha) ) < 3 \log{k}
\end{equation} 
holds for all $ k\geq 2 ,$ which will be useful in our study.

In \cite{16}, Bravo and Luca showed that
\begin{equation*}
L_n^{(k)} = \sum_{i=1}^k  (2 \alpha_{i}-1)   f_k (\alpha_i) (\alpha_{i})^{n-1} \quad \text{and} \quad \left \lvert L_n^{(k)} -  f_k (\alpha) (2 \alpha-1) \alpha^{n-1}   \right \rvert < 3/2,
\end{equation*}
for all $ k\geq 2. $ Thus,
\begin{equation}
\label{kLe}
L_n^{(k)} =  f_k (\alpha) (2 \alpha-1) \alpha^{n-1} + e_k(n),
\end{equation}
where $ \lvert e_k(n) \rvert <3/2. $ As in the classical case where $ k=2 $, we have the bounds
\begin{equation}
\label{L1}
\alpha^{n-1} \leq L_n^{(k)} \leq 2\alpha^{n}, 
\end{equation}
for all $n\geq 1$ and $ k \geq 2 $ \cite{16}. It is known that
\begin{equation}\label{nlessk}
L_i^{(k)} =3 \cdot 2^{i-2}, \quad \text{ for all } \quad 2 \leq i \leq k.
\end{equation}

\section{Proof of Theorem \ref{main1}}
First, we may directly derive some useful relations among the variables in \eqref{LnLmLlR}. Indeed, from the fact that 
$$ 10^{d-1} < a \left( \dfrac{10^{d}-1}{9} \right) < 10^d ,$$
together with \eqref{L1}, we find the following two relations
\begin{equation}\label{d1}
\left( \dfrac{1+\sqrt{5}}{2} \right)^{n-3} \leq \alpha ^{n+m+l-3} < 10^d  \Rightarrow \dfrac{n-3}{5}< d,
\end{equation}
and
\begin{equation}\label{d2}
10^{d-1} < 8 \alpha ^{n+m+l} < 2^{n+m+l+3}  \Rightarrow d < n+2 \leq 2n,
\end{equation}
for all $ n \geq 2. $

Before further calculations, by a computer search, we checked all repdigits having at least two digits and belonging to the products $ L_n^{(k)}L_m^{(k)}L_l^{(k)} $ in the range $ 2 \leq k \leq 25 $ and $ 0 \leq l \leq m \leq n \leq 25 ,$ to find the variables given in Theorem \ref{main1}. So from now on, we will take $ \max\{ k, n \} > 25.$

\subsection{The Case $ n \leq k. $ }
By substituting \eqref{nlessk} in \eqref{LnLmLlR}, we get
\begin{equation*}
3^5 \cdot 2^{n+m+l-6} = a ({10^{d}-1}).
\end{equation*}
Then, as $ {10^{d}-1} $ is odd, we find $ n+m+l-6  \leq 3, $ and hence we get $ n \leq 9 .$  Since $ L_9^{(k)} \leq 384 $ for all $ k \geq 2, $ by a computer search, we found that there is no solution of \eqref{LnLmLlR} when  $ n \leq k. $

Now assume that $ n \geq k+1 .$
\subsection{The Case $n \geq k+1 .$ }
Note that, since $ \max\{ k, n \}=n > 25,$ from \eqref{d1}, we may take $ d > 4. $ First, by \eqref{kLe}, we rewrite \eqref{LnLmLlR} as
\begin{multline*}
\left( f_k (\alpha) (2 \alpha-1) \alpha^{n-1} + e_k(n) \right)\left( f_k (\alpha) (2 \alpha-1) \alpha^{m-1} + e_k(m) \right) \\
\left( f_k (\alpha) (2 \alpha-1) \alpha^{l-1} + e_k(l) \right) = a10^{d}/9 - a/9.
\end{multline*}
Thus, we may write
\begin{equation*}
f_k (\alpha)^3 (2 \alpha-1)^3 \alpha^{n+m+l-3} - a10^{d}/9 = U_1-a/9,
\end{equation*}
where
\begin{align*}
U_1= & f_k (\alpha)^2 (2 \alpha-1)^2 \alpha^{-2} \left( \alpha^{n+m}e_k(l)+\alpha^{n+l}e_k(m)+\alpha^{m+l}e_k(n) \right) +\\
& f_k (\alpha) (2 \alpha-1) \alpha^{-1} \left( \alpha^{n}e_k(l)e_k(m)+\alpha^{m}e_k(n)e_k(l)+\alpha^{l}e_k(n)e_k(m) \right)+\\
& e_k(n)e_k(m)e_k(l).
\end{align*}
Then, by dividing both sides by $ f_k (\alpha)^3 (2 \alpha-1)^3 \alpha^{n+m+l-3} ,$ and taking the absolute values, we get
\begin{equation}
\label{Lam1}
\left \lvert \Lambda_1  \right \rvert <  \dfrac{15}{\alpha^{l-3}},
\end{equation}  
where
\begin{equation*}
\Lambda_1  := \alpha^{-(n+m+l-3)} 10^d  f_k (\alpha)^{-3} (2 \alpha-1)^{-3}    \dfrac{a}{9} -1. 
\end{equation*}
To get the above result, we used the facts that
\begin{equation*}
\dfrac{ \left \lvert U_1 \right \rvert}{f_k (\alpha)^3 (2 \alpha-1)^3 \alpha^{n+m+l-3}} \leq 3 \dfrac{3/2}{\alpha^{l-1}} +3 \dfrac{9/4}{\alpha^{2l-2}} + \dfrac{27/8}{\alpha^{n+m+l-3} } \leq \dfrac{107/8}{\alpha^{l-3}}
\end{equation*} 
and 
\begin{equation*}
\dfrac{a}{9}\dfrac{1}{ f_k (\alpha)^3 (2 \alpha-1)^3 \alpha^{n+m+l-3} }  \leq \dfrac{1}{\alpha^{l-3}}.
\end{equation*}
Also, we were taking into account that $ 1< f_k (\alpha) (2 \alpha-1) .$

We turn back to \eqref{LnLmLlR}, and write it as
\begin{equation*}
\left( f_k (\alpha) (2 \alpha-1) \alpha^{n-1} + e_k(n) \right)\left( f_k (\alpha) (2 \alpha-1) \alpha^{m-1} + e_k(m) \right) L_l^{(k)} = a10^{d}/9 - a/9.
\end{equation*}
Thus, we may write
\begin{equation*}
f_k (\alpha)^2 (2 \alpha-1)^2 \alpha^{n+m-2} L_l^{(k)} - a10^{d}/9 = U_{2}L_l^{(k)}-a/9,
\end{equation*}
where
\begin{equation*}
U_2=f_k (\alpha) (2 \alpha-1) \alpha^{-1} \left( \alpha^{n}e_k(m)+\alpha^{m}e_k(n) \right) +\\
e_k(n)e_k(m).
\end{equation*}
This time, we divide both sides by $ f_k (\alpha)^2 (2 \alpha-1)^2 \alpha^{n+m-2} L_l^{(k)} , $ and get
\begin{equation}
\label{Lam2}
\left \lvert \Lambda_2  \right \rvert <  \dfrac{25/4}{\alpha^{m-2}},
\end{equation}  
where
\begin{equation*}
\Lambda_2  := \alpha^{-(n+m-2)} 10^d  f_k (\alpha)^{-2} (2 \alpha-1)^{-2}    {L_l^{(k)}}^{-1} \dfrac{a}{9} -1. 
\end{equation*}
To get the above result, we used the facts that
\begin{equation*}
\dfrac{\lvert U_2 \rvert}{f_k (\alpha)^2 (2 \alpha-1)^2 \alpha^{n+m-2}} \leq   2 \dfrac{3/2}{\alpha^{m-1}} +\dfrac{9/4}{\alpha^{n+m-2}} \leq \dfrac{21/4}{\alpha^{m-2}}
\end{equation*}
and 
\begin{equation*}
\dfrac{a}{9}\dfrac{1}{ f_k (\alpha)^2 (2 \alpha-1)^2 \alpha^{n+m-2} L_l^{(k)} }  \leq \dfrac{1}{\alpha^{m-2}}.
\end{equation*}
We rearrange \eqref{LnLmLlR} again by rewriting it as
\begin{equation*}
\left( f_k (\alpha) (2 \alpha-1) \alpha^{n-1} + e_k(n) \right) L_m^{(k)}  L_l^{(k)} = a10^{d}/9 - a/9.
\end{equation*}
So, we have that
\begin{equation*}
f_k (\alpha) (2 \alpha-1) \alpha^{n-1} L_m^{(k)}  L_l^{(k)} - a10^{d}/9 = -  e_k(n) L_m^{(k)}  L_l^{(k)}  - a/9.
\end{equation*}
Dividing by $ f_k (\alpha) (2 \alpha-1) \alpha^{n-1} L_m^{(k)}  L_l^{(k)} $ and taking the absolute values, we get that
\begin{equation}
\label{Lam3}
\left \lvert \Lambda_3  \right \rvert <  \dfrac{5/2}{\alpha^{n-1}},
\end{equation}  
where
\begin{equation*}
\Lambda_3  := \alpha^{-(n-1)} 10^d  f_k (\alpha)^{-1} (2 \alpha-1)^{-1}   {L_m^{(k)}}^{-1} {L_l^{(k)}}^{-1} \dfrac{a}{9} -1. 
\end{equation*}
To get the above result, we used only the fact that
\begin{equation*}
\dfrac{ \lvert e_k(n) \rvert }{ f_k (\alpha) (2 \alpha-1) \alpha^{n-1} }+ \dfrac{a}{9} \dfrac{1}{f_k (\alpha) (2 \alpha-1) \alpha^{n-1} L_m^{(k)}  L_l^{(k)} } \leq \dfrac{3/2+1}{\alpha^{n-1}} = \dfrac{5/2}{\alpha^{n-1}}.
\end{equation*}
For each $ \Lambda_1, $ $ \Lambda_2 $ and $ \Lambda_3 ,$ we will apply Theorem \ref{Matveev} by taking $ t=3 $. In each case, we take $ \eta_1:=\alpha,$  $\eta_2:=10 $, $ b_2:=d ,$ $ b_3:=1 ,$   $ \mathbb{K}=\mathbb{Q}(\alpha) $ which is a real number field  with degree  $ d_\mathbb{K}=k .$ Also, $ ( \lvert b_1 \rvert,\eta_3)=( n+m+l-3,  f_k (\alpha)^{-3} (2 \alpha-1)^{-3} (a/9)  )  $ for $ \Lambda_1, $ $ ( \lvert b_1 \rvert,\eta_3)=( n+m-2,   f_k (\alpha)^{-2} (2 \alpha-1)^{-2}(a/9){L_l^{(k)}}^{-1}  )  $ for $ \Lambda_2 $ and $ ( \lvert b_1 \rvert,\eta_3)=( n-1, f_k (\alpha)^{-1} (2 \alpha-1)^{-1}  (a/9) {L_m^{(k)}}^{-1} {L_l^{(k)}}^{-1} )  $ for $ \Lambda_3. $ So, from \eqref{d2}, we will take $ B $ as $ 3n,2n $ and $ n+2 $ for $ \Lambda_1, $ $ \Lambda_2 $ and $ \Lambda_3 ,$ respectively. 

Since $h(\alpha )=({1}/{k})\log{\alpha } $ and  $  h(10)= \log {10},$  we take $ A_1 =  \log \alpha $ and $ A_2 = k \log{10}$ for each $ \Lambda_i. $ But the parameter $ A_3 $ will be different for each $ \Lambda_i $ and therefore we calculate it for each $ \Lambda_i $ separately.

To calculate the logarithmic height of $ \eta_3, $ we use the two facts, $ h( f_k(\alpha) ) < 3 \log{k}$ from \eqref{hfk}, and $ h( 2 \alpha -1 ) < \log{3}$ from \cite[page 147]{16}. Thus, we find
\begin{align*}
h((a/9) f_k (\alpha)^{-3} (2 \alpha-1)^{-3}) &  \leq  h(a/9) + 3 h( f_k(\alpha) ) + 3 h(2 \alpha -1 ) \\
& \leq  \log{9} + 9 \log{k} +  3 \log{3} \\
& \leq \log{243} + 9 \log{k} < 17 \log{k}.
\end{align*}
So, we take $ A_3:= 17k \log{k}$ for $ \Lambda_1. $ Now we will show that $ \Lambda_1 \neq 0. $ Assume that $ \Lambda_1 = 0. $ Then, we get
\begin{equation*}
 (a/9) 10^d =  f_k (\alpha)^{3} (2 \alpha-1)^{3} \alpha^{(n+m+l-3)}.   
\end{equation*}
Now, we take the image of both sides of this relation by applying any one of the automorphisms  $ \sigma_i : \alpha \rightarrow \alpha_i $ for any $ i \geq 2 .$ Then, by \eqref{fkprop}, we take the absolute values to get
$$
10^3 \leq (a/9) 10^d  \leq \lvert f_k(\alpha_i) \rvert ^3  \lvert  2 \alpha_i -1 \rvert^{3}  \lvert  \alpha_i  \rvert^{ (n+m+l-3) } \leq 27,
$$
which is clearly false. So, $ \Lambda_1 \neq 0 $ and hence we may apply Theorem \ref{Matveev} to \eqref{Lam1} which gives us a bound for $ \Lambda_1 .$ On the other hand, from \eqref{Lam1}, we directly see that 
\begin{equation*}
\log{ \lvert  \Lambda_1 \rvert } <  \log {15} - (l-3) \log {\alpha}.  
\end{equation*}
We combine these two results and get the following:
$$
(l-3) \log {\alpha} - \log {15} < -1.4 \cdot 30^6 \cdot 3^{4.5} \cdot k^2 (1+\log {k})(1+\log {3n}) \log \alpha \cdot k\log {10}    \cdot 17k \log {k}.
$$
Note that
$$ 1+\log {k} < (5/2)\log {k} \quad \text{and} \quad  1+\log {3n} < 1.6 \log {3n}   $$
hold for all $ k \geq 2 $ and $ n \geq 2. $ By taking into account these two facts, we obtain that
\begin{equation}
\label{l1}
l <  2.3 \cdot 10^{13} \cdot k^4 \log^2 {k} \log {3n}.
\end{equation}
For the next step, we need to know that
$ \Lambda_2 \neq 0. $ Assume $ \Lambda_2 = 0. $ Then, we get
\begin{equation*}
 (a/9) 10^d = f_k (\alpha)^{2} (2 \alpha-1)^{2} \alpha^{(n+m-2)} {L_l^{(k)}}.
\end{equation*}
Similarly for $ \Lambda_1 ,$  we get that
$$
(a/9) 10^d  \leq  \lvert f_k(\alpha_i) \rvert ^2  \lvert  2 \alpha_i -1 \rvert^{2} {L_l^{(k)}} \leq 9 {L_l^{(k)}}.
$$
This inequality together with \eqref{LnLmLlR}, implies that
\begin{equation*}
L_n^{(k)} L_m^{(k)}L_l^{(k)} \leq 9 {L_l^{(k)}} - (a/9) < 9 {L_l^{(k)}}. 
\end{equation*}
Thus, we find that $ L_n^{(k)} \leq 9, $ which is a contradiction to the choice of $ n>25. $ So $ \Lambda_2 \neq 0 .$ Since
\begin{align*}
h(  (a/9)f_k (\alpha)^{-2} (2 \alpha-1)^{-2}{L_l^{(k)}}^{-1}  ) & \leq  h(a/9) + 2 h( f_k(\alpha) ) + 2 h(2 \alpha -1 ) + h\left(  L_l^{(k)} \right) \\
& \leq  \log{9} + 6 \log{k} +  2 \log{3} + \log{2 \alpha^l}\\
& < \log{162} + 6 \log{k} + l \log{ \alpha } <  14 \log{k} + l \log{2},
\end{align*}
by using the upper bound of $ l $ given in \eqref{l1}, we find that
$$ h(  (a/9)f_k (\alpha)^{-2} (2 \alpha-1)^{-2}{L_l^{(k)}}^{-1}  ) \leq   1.6 \cdot 10^{13} \cdot k^4 \log^2 {k} \log {3n}.$$
Thus, we may take
$$ A_3 :=  1.6 \cdot 10^{13} \cdot k^5 \log^2 {k} \log {3n}  $$
for $ \Lambda_2. $ Now, we are ready to apply Theorem \ref{Matveev} to \eqref{Lam2}. By using the bound
\begin{equation*}
\log{ \lvert \Lambda_2 \rvert } <  \log {(25/4)} - (m-2) \log {\alpha}, 
\end{equation*}
we obtain
\begin{equation}
\label{m1}
m <  2.2 \cdot 10^{25} \cdot k^8 \log^3 {k} \log^2 {3n}.
\end{equation}
To find the value of $ A_3 $ for $ \Lambda_3, $ we need the following logarithmic height

\begin{align*}
h( (\dfrac{a}{9})f_k (\alpha)^{-1} (2 \alpha-1)^{-1}   {L_m^{(k)}}^{-1} {L_l^{(k)}}^{-1} ) & \leq  h(\dfrac{a}{9}) +  h( f_k(\alpha) ) +  h(2 \alpha -1 )\\
& + h\left(  L_l^{(k)}  L_m^{(k)}\right) \\
& \leq   \log{9} + 3 \log{k} +   \log{3} + \log{4 \alpha^{l+m} }\\   
& \leq \log{108} + 3 \log{k} + (l+m) \log{2} \\
& < 10 \log{k} + 2m \log{2}.
\end{align*}
Therefore, by using the bound of $ m $ obtained in \eqref{m1}, we take
$$ A_3 :=  3.1 \cdot 10^{25} \cdot k^9 \log^3 {k} \log {3n},$$
for $ \Lambda_3. $ The similar argument to the one that we used for $ \Lambda_2 ,$ shows that $ \Lambda_3 \neq 0 .$ Thus, by applying Theorem \ref{Matveev} to $ \Lambda_3 ,$ together with the inequality \eqref{Lam3}, we obtain that
\begin{equation*}
3n <  1.3 \cdot 10^{38} \cdot k^{12} \log^4 {k} \log^3 {3n}.
\end{equation*}
We need the following lemma from \cite{Guzman} to get a bound of $ n $ depending on $ k .$
\begin{lemma}
Let $ s \geq 1 $ and $ T>(4s^2)^s .$ If
$$
\dfrac{x}{( \log {x} )^s} < T, \quad \text{then} \quad x < 2^sT(\log{T})^s.
$$
\end{lemma}
We take $ s=3 $ and  $ T:= 1.3 \cdot 10^{38} \cdot k^{12} \log^4 {k} .$ Then,
\begin{align*}
(\log(T))^3  & < (\log(1.3) +38\log10 +12\log{k} +4\log( \log {k}) )^3 \\
& <  170^3 \log^3{k}.
\end{align*} 
Thus, we get the following  bound :
\begin{equation}\label{nfirst}
n< 1.8 \cdot 10^{45} k^{12} \log^7{k}.
\end{equation}
Now, we treat the cases $ k \leq 650 $ and $ k > 650 $ separately.
\subsection{The Case $ k \leq 650 $}

Let $ 2 \leq k \leq 650. $ Then, from \eqref{nfirst},  $ n $ is also bounded. Let
$$ \Gamma_1 := -(n+m+l-3) \log \alpha + d \log 10 + \log(f_k(\alpha)^{-3} (2 \alpha-1)^{-3} \cdot (a/9) ).$$
Then 
$$ \lvert \Lambda_1 \rvert := \left \lvert \exp(\Gamma_1) -1 \right \rvert < 15/ \alpha^{l-3}.$$
We claim that ${l} \leq 300.$ Suppose that this is not the case. Then, $ 15/ \alpha^{l-3}<1/2 $ and therefore $\lvert \Gamma_1 \rvert < 30/ \alpha^{l-3}.$  So, we have
\begin{equation}
\label{g1}
0 < \left \lvert (n+m+l-3) \dfrac{\log \alpha }{\log 10 }  - d + \mu_{(k,a)}   \right \rvert < 30/ \alpha^{l-3} \log 10.
\end{equation}
where
$$ 
\mu_{(k,a)}  :=  - \dfrac{ \log(f_k(\alpha)^{-3} (2 \alpha-1)^{-3} \cdot (a/9)) }{\log 10 } .
$$
For each $ 2 \leq k \leq 650 $, we take $ M_k := 5.4 \cdot 10^{45} k^{12} \log^7{k} >3n > n+m+l-3 $ and $\tau_k = \dfrac{ \log{\alpha} }{\log 10} .$ Also,  for each $ k $, we find a convergent $ p_i/q_i $ of the continued fraction of irrational number  $\tau_k,$ such that $ q_i > 6M_k.$ Then, we calculate
$$\epsilon_{(k,a)} := \lvert\lvert \mu_{(k,a)}  q_{i} \rvert\rvert -M_k \lvert\lvert \tau_k q_{i} \rvert\rvert  $$
for each $ a \in \{1, 2, \ldots , 9 \}.$ 
If $\epsilon_{(k,a)}< 0 ,$ then we repeat the same calculations for $ q_{i+1}. $ For each $ k ,$ except for $(k,a)=(2,9),$ we found an appropriate  $q_i,$ such that $ \epsilon_{(k,a)} > 0. $ In fact,  $ 0.000308  <\epsilon_{(k,a)}  .$

Thus, from  Lemma \ref{reduction}, we find an upper bound on $l-3 $  for each $ 2 \leq k \leq 650 $ such that none of these bounds is greater than $ 297. $ So, when $(k,a) \neq (2,9),$ we conclude that  $ l \leq 300 ,$ as we claimed.

If $(k,a) = (2,9),$ then by taking into account $ \alpha^2 =\alpha +1 $ for $ k=2, $  we see that $ 3 \log{\alpha} + \log(f_2(\alpha)^{-3} (2 \alpha-1)^{-3} =0  ,$ and hence $\epsilon_{(2,9)}= 0 $ for all $ q_i. $ Thus, Lemma \ref{reduction} can not be applicable in this case. Now, we will use directly the properties of continuous fractions. From \eqref{g1}, we get
\begin{equation}\label{cf1}
\left \lvert \dfrac{\log \alpha }{\log10} - \dfrac{d}{n+m+l}  \right \rvert < \dfrac{30}{ (n+m+l) \alpha^{l-3} \log 10 }  .
\end{equation}
If $ \dfrac{30}{ (n+m+l) \alpha^{l-3} \log 10 } <   \dfrac{1}{ 2(n+m+l)^2 } ,$ then $ \dfrac{d}{n+m+l} $ is a convergent of continued fraction expansion  of the irrational number $ \dfrac{\log \alpha }{\log10}  ,$ say $\dfrac{p_i}{q_i}.$ Since $ p_i $ and $ q_i $ are relatively prime, we deduce that $ q_i \leq n+m+l \leq 3n < 3 \cdot 5.7 \cdot 10^{47}< 1.8 \cdot 10^{48}. $ 

By Maple, we see  that $ i<101. $ Let $[a_0,a_1,a_2,a_3,a_4,\ldots]=[0,4,1,3,1,1,1,6, \ldots]$ be the continued fraction expansion of  $   \dfrac{\log \alpha }{\log10} .$ Then,  $ \max\{a_i\}  = 106$ for $ i=0,1,2, \ldots ,101. $ Thus, from the well-known property of continued fractions, see for example \cite[Theorem 1.1.(iv)]{Hen},  we write
$$  \dfrac{1}{106 \cdot (n+m+l)^2 } \leq    \dfrac{30}{ (n+m+l) \alpha^{l-3} \log 10 }   .$$
Therefore, from the inequality
$$ \alpha^{l-3} < \dfrac{106 \cdot 30 \cdot (n+m+l) }{\log 10} <  2.5 \cdot 10^{52}, $$ 
we find that $ l < 250.  $  This bound is valid also in the case
$$ \dfrac{30}{ (n+m+l) \alpha^{l-3} \log 10 } \geq   \dfrac{1}{ 2(n+m+l)^2 } .$$
So, whether $ (k,a)=(2,9) $ or not, we have that $ l \leq 300 .$

Let
$$ \Gamma_2 := -(n+m-2) \log \alpha + d \log 10 + \log(f_k(\alpha)^{-2} (2 \alpha-1)^{-2}  {L_l^{(k)}}^{-1} (a/9)) .$$
Then, 
$$ \lvert \Lambda_2 \rvert := \left \lvert \exp(\Gamma_2) -1 \right \rvert < (25/4)/\alpha^{m-2}.$$
Next, we show that there is a bound for $ m .$ Suppose that $ m>10. $
Then, $ (25/4)/\alpha^{m-2}<1/2 $ and therefore $\lvert \Gamma_2 \rvert < (25/2)/ \alpha^{m-2}.$  So, we have
\begin{equation}
\label{g2}
0 < \left \lvert (n+m-2) \dfrac{\log \alpha }{\log 10 }  - d +\mu_{(k,l,a)}  \right \rvert < 25/ ( 2 \alpha^{m-2} \log 10),
\end{equation}
where
$$ \mu_{(k,l,a)} = - \dfrac{ \log(f_k(\alpha)^{-2} (2 \alpha-1)^{-2} {L_l^{(k)}}^{-1} (a/9)) }{\log 10 }  .$$
For each $ 2 \leq k \leq 650 $, we take $ M_k := 3.6 \cdot 10^{45} k^{12} \log^7{k} > 2n > n+m-2 $ and $\tau_k = \dfrac{ \log{\alpha} }{\log 10} .$ Also, for each $ k ,$ like in the previous calculations,  we find an appropriate $ q_i > 6M_k.$ Then, we calculate
$\epsilon_{(k,l,a)} := \lvert\lvert \mu_{(k,l,a)}  q_{i} \rvert\rvert -M_k \lvert\lvert \tau_k q_{i} \rvert\rvert  $ for each $ 1 \leq a \leq 9 $ and  $ 0 \leq l \leq 300.$  We find that $ \epsilon_{(k,l,a)} > 0.132192 \cdot 10^{-6} $ for all $ 3 \leq k \leq 650.$

Thus, from  Lemma \ref{reduction}, we find an upper bound on $m-2 ,$  for each $ 3 \leq k \leq 650 ,$ and none of these bounds is greater than $ 298. $

Except for $(k,l,a)=(2,1,9),$ if $ k=2 ,$  then $ \epsilon_{(2,l,a)} > 0.327364 \cdot 10^{-38} ,$ and hence from  Lemma \ref{reduction}, $ m \leq 606. $

If $(k,l,a)=(2,1,9),$ then $ 2 \log{\alpha} + \log(f_2(\alpha)^{-2} (2 \alpha-1)^{-2} )=0 $ and hence, from \eqref{g2} we find that
\begin{equation*}
\left \lvert \dfrac{\log \alpha }{\log10} - \dfrac{d}{n+m}  \right \rvert < \dfrac{(25/2)}{ (n+m) \alpha^{m-2} \log 10 } .
\end{equation*}
Now, we follow the same argument as we did for \eqref{cf1} to show that $ m<250. $ So, we have that $ m \leq 300 $ if  $ 3 \leq k \leq 650,$ and  $ m \leq 606 $ if $ k=2. $

Let
$$ \Gamma_3 := -(n-1) \log \alpha + d \log 10 + \log(f_k(\alpha)^{-1} (2 \alpha-1)^{-1}  {L_l^{(k)} }^{-1} {L_m^{(k)} }^{-1} (a/9)) .$$
Thus,
$$ \lvert \Lambda_3 \rvert := \left \lvert \exp(\Gamma_3) -1 \right \rvert < 5/(2 \alpha^{n-1}) < (1/2).$$
Then $\lvert \Gamma_3 \rvert < 5/ \alpha^{n-1},$  and we have
\begin{equation}
\label{g3}
0 < \left \lvert (n-1) \dfrac{\log \alpha }{\log 10 }  - d     \right \rvert < 5/ (  \alpha^{n-1} \log 10),
\end{equation}
where
$$\mu_{(k,l,m,a)} = - \dfrac{ \log(f_k(\alpha)^{-1} (2 \alpha-1)^{-1}  {L_l^{(k)} }^{-1} {L_m^{(k)} }^{-1} (a/9))  }{\log 10 }. $$
This time, we take $ M_k := 1.8 \cdot 10^{45} k^{12} \log^7{k} > n-1  $ and $\tau_k = \dfrac{ \log{\alpha} }{\log 10} .$ For each $ k ,$ like in the previous calculations,  we find an appropriate $ q_i > 6M_k.$ Then, we calculate
$\epsilon_{(k,l,m,a)} := \lvert\lvert \mu_{(k,l,m,a)}  q_{i} \rvert\rvert -M_k \lvert\lvert \tau_k q_{i} \rvert\rvert  $ for each $ 1 \leq a \leq 9 $ and  $ 0 \leq l \leq m \leq 300 .$  This is done by taking into account the facts that
\begin{equation}
\label{ordp}
v_2\left( L_l^{(k)} L_m^{(k)} \right) \leq v_2(a) \quad \text{and} \quad  v_5\left( L_l^{(k)} L_m^{(k)} \right) \leq v_5(a)
\end{equation}
to reduce the unnecessary computations.
Except for $(k,l,m,a)=(2,1,1,9),$ we find  $ \varepsilon_k:=\min\{ \epsilon_{(k,l,m,a)} \} $ for each $ k .$ Thus, from  Lemma \ref{reduction}, we obtain an upper bound for $ n-1 $, and hence for $ n ,$ say $ n(k) ,$  for each $ 2 \leq k \leq 650 .$ Indeed, none of these bounds are greater than $ 297.5 $ when $ k \geq 3. $ Some of these bounds are $ n(3)<230.93 ,n(100)<263.5$ $ n(200)<273.95 ,$ $n(300)<284.7,$ $n(400)<281.1,$ $n(500)<284.9,$ $n(600)<293.97,$ $n(605)<297.4$ and $n(650)<290.3.$ Whereas for $ k=2 ,$ $n(2)<1216.$ In fact, even though $ q_{101}>6M_k:=6M_2,$ we have to take $ q_{305} $ for the condition $\epsilon_{(2,l,m,a)}>0$ and we find that  $ \epsilon_{(2,l,m,a)}> 0.134539 \cdot 10^{-91} .$ 

If $(k,l,m,a) = (2,1,1,9) ,$ then $  \log{\alpha} + \log(f_2(\alpha)^{-1} (2 \alpha-1)^{-1} )=0 $ and hence, from \eqref{g3}, we find that
\begin{equation*}
0 < \left \lvert \dfrac{\log \alpha }{\log10} - \dfrac{d}{n}  \right \rvert < \dfrac{5}{ n \alpha^{n-1} \log 10 } .
\end{equation*}
Now, we follow the same calculations as we did for \eqref{cf1} to show that the bound $ n(2)<1216$ is valid.

Finally, we write a short computer programme to check that the variables $ n,m,k,l,a $ and $ d $ are satisfying \eqref{LnLmLlR} by using the bounds $ n(k) $ for $ 2\leq k \leq 650 $ together with   \eqref{ordp} and \eqref{d2}. As a result, we find that there is no new solution of \eqref{LnLmLlR} except for those that were given in Theorem \ref{main1}.

\subsection{The Case $ k > 650 $ }

We cite the following lemma  from \cite[Lemma 2.6]{24}.
\begin{lemma}\label{Guzman2}
If $ k \geq 2 $ and $ n \geq k+1$ then 
$$  L_n^{(k)}=3 \cdot 2^{n-2}(1+\zeta(n,k)), \quad \text{where} \quad \lvert \zeta(n,k) \rvert  < \dfrac{1}{2^{k/2}} .$$
\end{lemma}
For $ k > 650 ,$ from \eqref{nfirst}, the inequality $ n< 2^{k/2} $ holds and hence from Lemma \ref{Guzman2}, we have 
\begin{equation*}
L_n^{(k)} L_m^{(k)}L_l^{(k)} =27 \cdot 2^{n+m+l-6}(1+\zeta(n,k))^3.
\end{equation*}
Thus, by \eqref{LnLmLlR}, we get
\begin{equation*}
(a/9)10^{d} - 27 \cdot 2^{n+m+l-6}  =   27 \cdot 2^{n+m+l-6}  ( 3 \zeta(n,k) +3 \zeta(n,k)^2 +\zeta(n,k)^3 )+ (a/9).
\end{equation*}
Hence,
\begin{equation*}
0 < \lvert \Lambda_4 \rvert := \left \lvert 2^{-(n+m+l-6)}10^d 3^{-5} a -1 \right \rvert < \dfrac{1}{2^{k/2-3}},
\end{equation*}
since $ ( 3 \zeta(n,k) +3 \zeta(n,k)^2 +\zeta(n,k)^3 )<7/2^{k/2}.$

Now, we apply Theorem \ref{Matveev} to $ \Lambda_4 $ by taking  $ (\eta_1, b_1) =( 2,{-(n+m+l-6)} ),$ $ (\eta_2, b_2) =( 10,d ),$  $ (\eta_3, b_3) =( 3^{-5} a, 1 ),$ $ \mathbb{K}=\mathbb{Q} $ and $ B:=3n .$ Thus, we find that
\begin{equation*}
(\dfrac{k}{2}-3) \log {2} < 1.4 \cdot 30^6 \cdot 3^{4.5} \cdot (1+\log{3n}) \log {2} \cdot \log {10 } \cdot 5\log {3}.
\end{equation*}
From \eqref{nfirst}, we may take
\begin{align*}
1+\log {3n}  & < 1+\log {5.4} +23 \log {10^{2}} + 12 \log {k} +  7 \log {\log {k}} \\
& <  45  \log {k},
\end{align*}
and hence we get that $ k < 1.63 \cdot 10^{14} \log {k} ,$ which implies that
\begin{equation*} 
k < 6 \cdot 10^{15}. 
\end{equation*}
So, by \eqref{nfirst}, we have
\begin{equation*}
n<3.3 \cdot 10^{245}.
\end{equation*} 

\subsection{Reducing the Bound on k}
Let 
\begin{equation}\label{gL27}
\Gamma_{4} := -(n+m+l-6) \log 2 + d \log 10 + \log(a/{3^5}).
\end{equation}
Then $ \vert \Lambda_{4} \rvert :=\left \lvert \exp(\Gamma_{4})-1 \right \rvert < \dfrac{1}{ 2^{k/2-3} } <\dfrac{1}{2}.$ Hence, we have that $ \lvert \Gamma_{4} \rvert < \dfrac{2}{2^{k/2-3}} .$ From \eqref{gL27}, we write
\begin{equation*}
0< \left \lvert (n+m+l-6) \dfrac{ \log 2 }{\log 10}  - d - \dfrac{ \log(a/{3^5}) }{ \log 10} \right \rvert < \dfrac{1}{ 2^{k/2-4}  \log {10} }.
\end{equation*}

Let $ M:=  10^{246} >3n>n+m+l-6 ,$ $\tau = \dfrac{ \log{2}}{\log 10} $ and $ \mu_a := -\dfrac{ \log(a/{3^5}) }{ \log 10}.$

Then, the denominator of the $ 504th $ convergent of $ \tau,$ say $q_{504},$ exceeds $ 6M$ and
$$ \epsilon_a := \lvert\lvert\mu_a q_{505} \rvert\rvert -M \lvert\lvert \tau q_{505} \rvert\rvert > 0.05055 ,$$
for each $ a \in \{1, 2, \ldots, 9 \}. $ Thus, by applying  Lemma \ref{reduction} with the parameters $ A:=\dfrac{1}{\log{10}} ,$ $ B:=2 $ and $ w:={k/2-4} $, we find that $ {k/2-4} < 825 .$ That is $ k <1660,$ and thus, by \eqref{nfirst}, $ n<9.8 \cdot 10^{89} .$ 

We repeat the same reduction steps by taking $ M:= 2.94 \cdot 10^{90} >3n > n+m+l-6 .$ Then, $q_{184}>6M $ and we find that
$$ \epsilon_a := \lvert\lvert\mu_a q_{184} \rvert\rvert -M \lvert\lvert \tau q_{184} \rvert\rvert >0.009382 .$$
By  Lemma \ref{reduction}, we find that $ {k/2-4} < 310 .$  That is $ k <630,$  which contradicts the fact $ k > 650. $ So, we conclude that Equation \eqref{LnLmLlR} has no solution when $ k > 650 .$ This completes the proof.



\end{document}